\documentclass[11pt]{amsart}
\usepackage{fullpage}
\usepackage{amsmath}
\usepackage{amsthm}
\usepackage{amssymb}
\usepackage{amsfonts,mathrsfs}
\usepackage{amsxtra}

\makeatletter
\def\@settitle{\begin{center}%
  \baselineskip14\p@\relax
  \bfseries
  \uppercasenonmath\@title
  \@title
  \ifx\@subtitle\@empty\else
     \\[1ex]\uppercasenonmath\@subtitle
     \footnotesize\mdseries\@subtitle
  \fi
  \end{center}%
}
\def\subtitle#1{\gdef\@subtitle{#1}}
\def\@subtitle{}
\makeatother

\usepackage{epsfig}
\usepackage{epstopdf}

\usepackage{url}
\usepackage[all]{xy}
\usepackage{latexsym}
\usepackage{amssymb}
\usepackage[cp850]{inputenc}
\usepackage{amsfonts}
\usepackage{graphicx}
\usepackage{dsfont}

\usepackage[usenames]{color}

\newtheorem{theorem}{Theorem}[section]

\newtheorem{prop}[theorem]{Proposition}

\newtheorem*{defi*}{Definition}			
\newtheorem*{bei*}{Example}
\newtheorem*{sat*}{Theorem}				
\newtheorem*{kor*}{Corollary}
\newtheorem*{rmk*}{Remark}				
	
\newtheorem*{quest*}{Question}

\let\ssection=\section
\renewcommand{\section}{\setcounter{equation}{0}\ssection}

\newtheorem*{namedtheorem}{\theoremname}
\newcommand{\theoremname}{testing}

\theoremstyle{remark}
\newtheorem*{bem}{Remark}

\newtheorem*{namedtheoremr}{\theoremnamer}
\newcommand{\theoremnamer}{testing}

\newcommand{\BR}{\mathbb R}			
			
			\newcommand{\BZ}{\mathbb Z}

\newcommand{\CC}{\mathcal C}		\newcommand{\calD}{\mathcal D}

		\newcommand{\CL}{\mathcal L}
\newcommand{\CM}{\mathcal M}

\newcommand{\FM}{\mathfrak m}

\newcommand{\actson}{\curvearrowright}
\newcommand{\D}{\partial}

\DeclareMathOperator{\Id}{Id}		

\DeclareMathOperator{\Map}{Map}

\newcommand{\comment}[1]{}

\DeclareMathOperator{\Stab}{Stab}

\DeclareMathOperator{\Homeo}{Homeo}

\DeclareMathOperator{\Thu}{Thu}

\begin{document}

\title{Mirzakhani's curve counting: from simple to all}
\author{Viveka Erlandsson}
\thanks{The first author gratefully acknowledges support from EPSRC grant EP/T015926/1}
\address{School of Mathematics, University of Bristol, Bristol BS8 1UG, UK {\rm and}\newline ${ }$ \hspace{0.2cm} Department of Mathematics and Statistics, UiT The Artic University of Norway}
\email{v.erlandsson@bristol.ac.uk}
\author{Juan Souto}
\address{IRMAR, Universit\'e de Rennes 1, Rennes, France}
\email{jsoutoc@gmail.com}
\begin{abstract}
In \cite{Maryam1}, Mirzakhani obtained the asymptotic growth, when $L\to\infty$, of the number of curves in the mapping class group orbit of some given simple curve and with length at most $L$. Years later she extended in \cite{Maryam2} this result from simple to arbitrary curves. Here we give a short and relative low-tech argument showing how to derive the general result from the one for simple curves.
\end{abstract}
\maketitle

\section{Counting curves}

Let $S$ be an orientable surface of negative Euler characteristics and of genus $g$, with $k$ ends, and with $\D S=\emptyset$, and suppose for the sake of concreteness that $S$ is not a thrice punctured sphere. We will also ignore the difference between essential closed curves---that is homotopically essential curves which are not homotopic into an end---and their homotopy classes. Under a {\em multicurve} we understand a finite formal sum $\gamma=\sum a_i\gamma_i$ of pairwise non-homotopic, primitive closed curves $\gamma_i$, with coefficients $a_i> 0$. The length $\ell_X(\gamma)$ of $\gamma$ with respect to a hyperbolic structure $X$ on $S$ is the weighted sum $\sum_i a_i\ell_X(\gamma_i)$ of the lengths of the associated geodesics. The multicurve $\gamma$ is {\em simple} if all the curves $\gamma_i$ are simple and disjoint. The mapping class group 
$$\Map(S)=\Homeo_+(S)/\Homeo_0(S)$$
of $S$ acts on the set of (homotopy classes of) curves, and hence on the set of multicurves. Two multicurves $\gamma$ and $\gamma'$ are {\em of the same type} if they belong to the same mapping class group orbit, meaning that there exists $\phi\in\Map(S)$ with $\gamma'=\phi(\gamma)$.

In this note we are interested in two results of Mirzakhani describing the asymptotic behavior of the number of multicurves of given type and length bounded by $L$, where length is measured with respect to, for instance, a hyperbolic metric on $S$. In \cite{Maryam1}, Mirzakhani proved the following result:

\begin{theorem}[Mirzakhani]\label{Mirzakhani-simple}
For any essential simple multicurve $\gamma$ in $S$, and for any complete, finite volume, hyperbolic structure $X$ on $S$ we have 
$$\lim_{L\to\infty}\frac{\#\{\gamma'\text{ of type }\gamma\text{ and with }\ell_X(\gamma')\le L\}}{L^{6g-6+2k}}=\frac{B(X)\cdot c(\gamma)}{m_{g,k}},$$
where $B(X)=\FM_{\Thu}(\{\lambda\in\CM\CL\text{ with }\ell_X(\lambda)\le 1\})$ is the Thurston measure of the set of measured laminations $\lambda\in\CM\CL(S)$ with $\ell_X(\lambda)\le 1$, where $m_{g,k}=\int B(X)$ is the Weil-Petersson integral of $B(X)$ over the moduli space, and where $c(\gamma)$ is a rational only depending on $\gamma$.
\end{theorem}

The constant $c(\gamma)$ in Theorem \ref{Mirzakhani-simple} is related to the Weil-Petersson volume of the set of hyperbolic surfaces where $\gamma$ has given length. Mirzakhani gives a recursive formula to compute these volumes using a generalization \cite{Maryam-McS} of McShane's identity \cite{McShane}, and it is from there that she derives that $c(\gamma)$ is rational. Simplicity of $\gamma$ plays a key r\^ole in these arguments. However, Mirzakhani later came up with another argument \cite{Maryam2} that allowed her to treat the non-simple case as well:

\begin{theorem}[Mirzakhani]\label{Mirzakhani-general}
For any essential multicurve $\gamma$ in $S$, and for any complete, finite volume, hyperbolic structure $X$ on $S$ we have
$$\lim_{L\to\infty}\frac{\#\{\gamma'\text{ of type }\gamma\text{ and with }\ell_X(\gamma')\le L\}}{L^{6g-6+2k}}=\frac{B(X)\cdot c(\gamma)}{m_{g,k}}$$
where $B(X)$ and $m_{g,k}$ are as in Theorem \ref{Mirzakhani-simple}, and where $c(\gamma)$ is again rational.
\end{theorem}

Besides Mirzakhani's original papers \cite{Maryam1,Maryam2} there are by now a few other sources for either all or  part of the material covered by these two theorems, and this is good news because, all our admiration and appreciation for Mirzakhani's work notwithstanding, the paper \cite{Maryam2} is very hard to read, containing some rather arcane parts (it has been suggested that Mirzakhani's argument in \cite{Maryam2} is only tight if the curve is filling---a bit of a moot point, if you ask us). Anyways, a proof of Theorem \ref{Mirzakhani-general} (and hence of Theorem \ref{Mirzakhani-simple}) appears in our monograph \cite{book}. The proof we present there has the twofold advantage that it applies directly to many other notions of length than the hyperbolic length, and that the argument is pretty coarse, meaning that only rather elementary tools are needed. The results and the arguments in \cite{EMM} and \cite{Francisco-III} are much more delicate. Indeed, in the first of those papers the authors give a different argument for Theorem \ref{Mirzakhani-simple}, an argument that not only yields the asymptotics but also an error term. Along the same lines, a proof for filling curves of Theorem \ref{Mirzakhani-general} with error term is given in \cite{Francisco-III}. 
\medskip

\underline{Our goal here is to explain how to derive Theorem \ref{Mirzakhani-general} from Theorem \ref{Mirzakhani-simple}}. It seems that doing so might be of some use because Mirzakhani's proof of Theorem \ref{Mirzakhani-simple} has over time been widely publicized and studied, and hence it makes sense to give an easy path from this theorem to the general Theorem \ref{Mirzakhani-general}. However, since detailed arguments are available elsewhere, it seems reasonable to make some additional assumptions allowing us to avoid noise that would be hiding the structure of the argument. 

\begin{quote}
{\bf Simplifying assumptions.} We assume that $S$ is a closed surface of genus $g\ge 3$. We will also just prove that the limit in Theorem \ref{Mirzakhani-general} exists (the reader would have no problem distilling out of the argument an expression for the limit). 
\end{quote}

In any case, in the final section of the paper we comment briefly on what has to additionally be done to deal with surfaces with cusps and/or smaller complexity.  
\medskip

For what it is worth, which is surely nothing, we are dedicating this paper to the people of Ukraine. We are horrified by what is being done to them.

\section{Geodesic currents}

Measured laminations were introduced by Thurston in the 70's and there are quite a few good sources for this material---we recommend \cite{Casson-Bleiler,Hatcher}. Currents were introduced by Bonahon in \cite{Bonahon1,Bonahon2,Bonahon3} and most of the facts that we will need can be found, in a more or less transparent way, in these three papers. In the case of closed surfaces, the one we care about here, \cite{Aramayona-Leininger} is a very readable account of currents, measured laminations, and the relation between them, but it is maybe understandable that we are partial for the presentation of these matters in our own book \cite{book}.
\medskip

Let $S$ be a closed hyperbolic surface of genus $g\ge 3$. {\em Geodesic currents}, or just plainly {\em currents}, on $S$ are $\pi_1(S)$-invariant Radon measures on the space of geodesics on the universal cover $\widetilde S$ of $S$. The space $\CC(S)$ of all currents on $S$ is endowed with the weak-*-topology. The reader might be surprised, but we will not really care about what curents are. Indeed, we will just care about some of the properties of the space $\CC(S)$ of currents. However, before coming to the list of properties of $\CC(S)$ recall that all curves are assumed to be essential, that a multicurve is a formal positive linear combination of curves (positive just means that the weights are positive), and that a multicurve is filling if its geodesic representative cuts the (closed) surface $S$ into disks. Here are the properties of $\CC(S)$:
\begin{enumerate}
    \item $\CC(S)$ is a locally compact metrizable topological space.
    \item $\CC(S)$ is a cone as a topological vector space, meaning in
      particular that there are continuous maps
    \begin{align*}
      \CC(S)\times\CC(S)&\to\CC(S),\ (\lambda,\mu)\mapsto\lambda+\mu\\
      \BR_{\ge 0}\times\CC(S)&\to\CC(S),\ (t,\lambda)\mapsto t\lambda
    \end{align*}
      satisfying the usual associativity, commutativity and distributivity
      properties as in vector spaces.
    \item The set $\{\gamma\text{ closed geodesic in }S\}$ is a subset
      of $\CC(S)$ and in fact the set \[\BR_+\cdot\{\gamma\text{
        closed geodesic in }S\}\] of weighted closed geodesics is
      dense in $\CC(S)$.
    \item The inclusion of the set of weighted simple geodesics into
      $\CC(S)$ extends to a continuous embedding of the space
      $\CM\CL(S)$ of measured laminations into $\CC(S)$.
    \item There is a continuous bilinear map
      $$\iota:\CC(S)\times\CC(S)\to\BR_{\ge 0}$$ such that
      $\iota(\gamma,\gamma')$ is nothing other than the geometric
      intersection number for all closed geodesics $\gamma,\gamma'$.
\item A multicurve $\gamma$ is filling if and only if $\iota(\gamma,\lambda)>0$ for all non-zero currents $\lambda$.
    \item The set $\{\lambda\in\CC(S)\text{ with
      }\iota(\lambda,\eta)\le L\}$ is compact for every $L\ge 0$ and every
      filling multicurve $\eta$. In particular, the projective space $P\CC(S)=(\CC(S)\setminus\{0\})/\BR_{>0}$ is compact. 
    \item The mapping class group acts continuously on $\CC(S)$ by
      linear automorphisms. Moreover, the inclusion of the set of closed
      geodesics into $\CC(S)$ is equivariant with respect to this
      action.
\item The action of the mapping class group on the open set of filling geodesic currents is properly discontinuous. 
  \end{enumerate}
As we already mentioned above, all of this is basically due to Bonahon, with the only exception being the last item in our list---this was proved in \cite{EM}. 
\medskip

Before moving on note that the fact that $\Map(S)$ acts properly discontinuously on the space of filling currents implies that for any current $\alpha\in\CC(S)$ the group $\Stab(\alpha)=\{\phi\in\Map(S)\text{ with }\phi(\alpha)=\alpha\}$ acts properly discontinuous on the set
$$\CC_\alpha=\{\lambda\in\CC(S)\text{ with }\lambda+\alpha\text{ filling}\}$$
because the action $\Stab(\alpha)\actson\CC_\alpha$ is conjugated to the action of $\Stab(\alpha)$ on the set $\{\lambda+\alpha\text{ with }\lambda\in\CC_\alpha\}$, a subset of the set of filling currents. If $\sigma$ is a filling multicurve in $S$ such that
$$\Stab(\sigma)=\{\phi\in\Map(S)\text{ with }\phi(\sigma)=\sigma\}=\Id$$
then the set
\begin{equation}\label{eq fundamental domain}
\calD_\alpha^\sigma=\{\lambda\in\CC_{\alpha}\text{ with }\iota(\lambda,\sigma)<\iota(\lambda,\phi(\sigma))\text{ for all }\phi\in\Stab(\alpha)\setminus\{\Id\}\}
\end{equation}
is an (open) fundamental domain for this action in the sense that 
\begin{itemize}
\item for every $\lambda\in\CC_{\alpha}$ there is at most one $\phi\in\Stab(\alpha)$ with $\phi(\lambda)\in\calD_\alpha^\sigma$, and
\item for every $\lambda\in\CC_{\alpha}$ there is at least one $\phi\in\Stab(\alpha)$ with $\phi(\lambda)$ in the closure
\begin{equation}\label{eq fundamental domain closure}
\bar{\calD}_\alpha^\sigma=\{\lambda\in\CC_{\alpha}\text{ with }\iota(\lambda,\sigma)\le\iota(\lambda,\phi(\sigma))\text{ for all }\phi\in\Stab(\alpha)\}
\end{equation} 
of $\calD_\alpha^\sigma$.
\end{itemize}
See \cite[Prop. 8.4]{book} for details.

\section{Convergence of measures}

Besides the spaces of measured laminations and currents, the other main player in our discussion is the Thurston measure, which we now introduce. First we need a bit of notation. Let $\CM\CL_\BZ(S)\subset\CM\CL(S)\subset\CC(S)$ be the set of all integrally weighted simple multicurves on $S$. Also, for $\alpha\in\CC(S)$, let $\delta_\alpha$ be the Dirac probability measure on $\CC(S)$ centered at $\alpha$. Now, the {\em Thurston measure} $\FM_{\Thu}$ is the Radon measure on  $\CC(S)$ given by 
$$\FM_{\Thu}=\lim_{L\to\infty}\frac 1{L^{6g-6}}\sum_{\alpha\in\CM\CL_\BZ(S)}\delta_{\frac 1L\alpha},$$
where $g$ still denotes the genus of our closed surface $S$, and where the convergence takes place with respect to the weak-*-topology on the space of measures on $\CC(S)$. More precisely, for any compact subset $K\subset\CC(S)$ with $\FM_{\Thu}(\D K)=0$ we have
$$\FM_{\Thu}(K)=\lim_{L\to\infty}\frac 1{L^{6g-6}}\#\left\{\alpha\in\CM\CL_\BZ(S)\text{ with }\frac 1L\alpha\in K\right\}.$$

\begin{bem}
The Thurston measure was originally described (up to a constant) as the measure associated to the Thurston's symplectic form on $\CM\CL(S)$---note that $\FM_{\Thu}$ is supported by $\CM\CL(S)$. As far as we know, the definition given here is due to Mirzakhani \cite{Maryam1}---see \cite{book} for a proof of the fact that the limit in the definition of $\FM_{\Thu}$ does actually exist and see \cite{Monin} for the relation between the symplectic description and the one above.
\end{bem}

In \cite{Maryam1} Mirzakhani considered, for a simple curve $\gamma$ the measure 
\begin{equation}\label{eq ynrohmai}
\FM_\gamma^L=\frac 1{L^{6g-6}}\sum_{\gamma'\text{ of type }\gamma}\delta_{\frac 1L\gamma'}.
\end{equation}
and studied its behavior as $L\to\infty$. 
Noting that 
$$\FM_\gamma^L\le\frac 1{L^{6g-6}}\sum_{\alpha\in\CM\CL_\BZ(S)}\delta_{\frac 1L\alpha}$$
and that the right hand side converges, we get that the measures $\{\FM_{\gamma}^L\}$ 
form a pre-compact family, meaning that any sequence $L_n\to\infty$ has a subsequence such that the limit $\lim_{i\to\infty}\FM_\gamma^{L_{n_i}}$ exists. 
Moreover, since the limit of the right hand side is the Thurston measure, we get that any measure $\FM$ of the form 
$$\FM=\lim_{n\to\infty}\FM_\gamma^{L_{n}}$$ 
is bounded by $\FM_{\Thu}$ and hence absolutely continuous with respect to it. This actually means that $\FM$ is a multiple of the Thurston measure: both $\FM$ and $\FM_{\Thu}$ are mapping class group invariant with $\FM_{\Thu}$ actually ergodic by Masur's ergodicity theorem \cite{Masur}. 

In this discussion we assumed that $\gamma$ was a simple curve. Passing from curves to multicurves involves just a little argument, but simplicty seems much more essential. This is probably one of the reasons why Mirzakhani approached in \cite{Maryam2} the general case differently, without considering the measures \eqref{eq ynrohmai}. However, in \cite{ourpaper, book} we proved, keeping the general philosophy of the argument above, that it all holds also for general multicurves: 

\begin{prop}\label{prop1}
Let $\gamma$ be a multicurve on $S$, simple or not. Any sequence $(L_n)_n$ of positive numbers with $L_n\to\infty$ has a subsequence $(L_{n_i})_i$ such that the measures $(\FM_{\gamma}^{L_{n_i}})_i$ converge in the weak-*-topology to the measure $c\cdot\FM_{\Thu}$ on $\CM\CL\subset\CC$ for some $c\ge 0$. 
\end{prop}

See \cite[Prop. 4.1]{ourpaper} or \cite{book} for a proof of Proposition \ref{prop1}. 
\medskip

The reason why we care about the measures $\FM_\gamma^L$ and their limits is the following:

\begin{prop}\label{prop2}
Let $\gamma$ be a possibly non-simple multicurve on $S$ and let $f_0:\CC(S)\to\BR_+$ be a positive, homogenous and continuous function. Then the following are equivalent:
\begin{enumerate}
\item The limit $\lim_{L\to\infty}\FM_\gamma^L$ exists,
\item the limit $\lim_{L\to\infty}\frac 1{L^{6g-6}}\#\{\gamma'\text{ of type }\gamma\text{ with }f(\gamma)\le L\}$ exists for every positive, homogenous and continuous function $f:\CC(S)\to\BR_{\ge 0}$, and
\item the limit $\lim_{L\to\infty}\frac 1{L^{6g-6}}\#\{\gamma'\text{ of type }\gamma\text{ with }f_0(\gamma)\le L\}$ exists.
\end{enumerate}
\end{prop}

Here a function $f:\CC(S)\to\BR_+$ is homogenous if we have $f(t\cdot\lambda)=t\cdot f(\lambda)$ for all $t\ge 0$ and $\lambda\in\CC(S)$ and positive if $f(\lambda)>0$ for all non-zero currents $\lambda$. There are plenty of examples of positive, homogenous and continuous functions $\CC(S)\to\BR_+$. For example, if $\sigma\in\CC(S)$ is a filling current then $\lambda\mapsto\iota(\lambda,\sigma)$ is such a function. Also, the hyperbolic length function on the set of curves extends to a positive, homogenous and continuous function $\ell_X:\CC(S)\to\BR_+$. These are the main examples we will care about in this note, but the reader can find many more such functions in \cite{hugo,didac-Thurston}.

\begin{proof}
Suppose that (1) holds. Then, from Proposition \ref{prop1} we get that there is $c\ge 0$ with 
\begin{equation}\label{eq prop2 eq1}
\lim_{L\to\infty}\FM_\gamma^L=c\cdot\FM_{\Thu}
\end{equation}
Now, let $f:\CC(S)\to\BR_{\ge 0}$ be a positive, homogenous and continuous function. Homogeneity implies that
\begin{align*}
\FM_\gamma^{L_n}(\{\alpha\in\CC(S)\text{ with }f(\alpha)\le 1\})
&=\frac 1{L^{6g-6}}\#\left\{\gamma'\text{ of type }\gamma\text{ with }f\left(\frac 1L\gamma\right)\le 1\right\}\\
&=\frac 1{L^{6g-6}}\#\left\{\gamma'\text{ of type }\gamma\text{ with }f\left(\gamma\right)\le L\right\}
\end{align*}
and positivity (and again homogeneity, and continuity, and the fact that $\CC(S)$ is locally compact) implies that these quantities are finite. Now, continuity and homogeneity of $f$, together with the fact that the Thurston measure is locally finite and satisfies $\FM_{\Thu}(tU) = t^{6g-6}\FM_{\Thu}(U)$ for all $U$, imply that $\FM_{\Thu}(\{f(\cdot)=1\})=0$. Armed with this we get from \eqref{eq prop2 eq1} that
$$\lim_{L\to\infty}\FM_\gamma^{L}(\{\alpha\in\CC(S)\text{ with }f(\alpha)\le 1\})=
c\cdot\FM_{\Thu}(\{\alpha\in\CC(S)\text{ with }f(\alpha)\le 1\})$$
Combining the two last equations we have proved that 
$$\lim_{L\to\infty}\frac 1{L^{6g-6}}\#\{\gamma'\text{ of type }\gamma\text{ with }f(\gamma)\le L\}=
c\cdot\FM_{\Thu}(\{\alpha\in\CC(S)\text{ with }f(\alpha)\le 1\})$$
and hence that (2) holds.

If (2) holds, then evidently (3) holds. It remains to prove that (3) implies (1). Noting that a precompact sequence converges if and only if every convergent subsequence has the same limit, and noting that from Proposition \ref{prop1} we get that the measures $\FM_\gamma^L$ form a precompact family when $L\to\infty$ and that every limit is of the form $c\cdot\FM_{\Thu}$, we deduce that to prove (1) it suffices to show that whenever we have
$$\lim_{n\to\infty}\FM_\gamma^{L_n}=c\cdot\FM_{\Thu}$$
then the constant $c$ does not depend on the sequence $(L_n)$. Well, applying the same computations as above, but this time to the sequence $(L_n)$ and to our given function $f_0$ we get that 
$$\lim_{n\to\infty}\frac 1{L_n^{6g-6}}\#\{\gamma'\text{ of type }\gamma\text{ with }f_0(\gamma)\le L_n\}=
c\cdot\FM_{\Thu}(\{\alpha\in\CC(S)\text{ with }f_0(\alpha)\le 1\})$$
On the other hand, the assumption (3) means that the limit on the left does not depend on the sequence $(L_n)$. It follows that 
$$c=\frac{\lim_{L\to\infty}\frac 1{L^{6g-6}}\#\{\gamma'\text{ of type }\gamma\text{ with }f_0(\gamma)\le L\}}{\FM_{\Thu}(\{\alpha\in\CC(S)\text{ with }f_0(\alpha)\le 1\})}$$
and this means in particular that $c$ does not depend on $(L_n)$, as we needed to prove. We are done.
\end{proof}

\section{From simple to all}\label{sec proof}
In this section we come to the goal of this note: {\underline{deduce Theorem \ref{Mirzakhani-general} from Theorem \ref{Mirzakhani-simple}}}.
\medskip

Suppose now that $\sigma$ is a filling multicurve and that it is generic in the sense that for all $\phi\in\Map(S)\setminus\{\Id\}$ we have
\begin{equation}\label{eq000}
\FM_{\Thu}(\{\lambda\in\CM\CL\text{ with }\iota(\lambda,\sigma)=\iota(\phi(\lambda),\sigma)\})=0.
\end{equation}
It is not hard to construct generic multicurves. For example, if $\gamma_1,\dots,\gamma_s$ are simple curves separating measured laminations in the sense that for $\lambda,\mu\in\CM\CL(S)$ we have
$$\iota(\lambda,\gamma_i)=\iota(\mu,\gamma_i)\text{ for all }i=1,\dots,s\ \Longrightarrow\ \lambda=\mu,$$ 
then we can take $\sigma=\sum a_i\gamma_i$ for generic vectors $(a_1,\dots,a_s)\in\BR_{>0}^s$---see \cite[Lemma 8.3]{book}. Anyways, note that genericity implies that the stabilizer $\Stab(\sigma)\subset\Map(S)$ of $\sigma$ in the mapping class group is trivial,
\begin{equation}\label{eq chichu is drinking}
\Stab(\sigma)=\Id
\end{equation}
 and hence that we can identify the mapping class group with the $\Map(S)$-orbit of $\sigma$.

As a first step in the proof we will deduce from Theorem \ref{Mirzakhani-simple} that the measures $\FM_\sigma^L$ converge when $L\to\infty$. As in the proof of Proposition \ref{prop2} it suffices to check that whenever 
\begin{equation}\label{eq1}
\lim_{n\to\infty}\FM_\sigma^{L_n}=c\cdot\FM_{\Thu}
\end{equation}
then $c$ does not depend on the particular sequence $(L_n)$. 

Well, to compute $c$ let us fix a simple multicurve $\alpha$ and note that from Theorem \ref{Mirzakhani-simple}, combined if you want with Proposition \ref{prop2}, we get that the limit
\begin{equation}\label{eq111}
\Lambda\stackrel{\text{def}}=\lim_{L\to\infty}\frac 1{L^{6g-6}}\vert\{\alpha'\text{ of type }\alpha\text{ with }\iota(\alpha',\sigma)\le L\}\vert
\end{equation}
exists. The name of the game now is to express $c$ in terms of $\Lambda$ and other quantities that might depend on $\alpha$ and $\sigma$ but do not depend on the particular sequence $(L_n)$. 

To do that let $\calD_\alpha^\sigma$, as in \eqref{eq fundamental domain}, be the fundamental domain for the action of $\Stab(\alpha)$ on 
$$\CC_\alpha=\{\lambda\in\CC\text{ such that }\lambda+\alpha\text{ is filling}\}$$
and as in and \eqref{eq fundamental domain closure} let $\bar\calD_\alpha^\sigma$ be the closure of $\calD_\alpha^\sigma$---this makese sense because our chosen $\sigma$ satisfies \eqref{eq chichu is drinking}. The reason why we care about this fundamental domain is that, by the mere act of being a fundamental domain, there is a set $\Theta\subset\Map(S)$ of representatives of the classes $\Map(S)/\Stab(\alpha)$ such that 
\begin{equation}\label{eq2}
\{\phi\in\Map(S)\text{ with }\phi^{-1}(\sigma)\in\calD_\alpha^\sigma\}\subset\Theta\subset\{\phi\in\Map(S)\text{ with }\phi^{-1}(\sigma)\in\bar\calD_\alpha^\sigma\}.
\end{equation}
The first inclusion in \eqref{eq2} yields:
\begin{align*}
\#\{ \alpha'\text{ of type }\alpha\text{ with }\iota(\alpha',\sigma)\le L_n\}
&=\#\left\{ \phi(\alpha)\text{ where }\phi\in\Map(S)\text{ and }\iota(\phi(\alpha), \sigma) \leq L_n\right\}\\
&=\#\left\{\phi\in\Theta\text{ with } \iota(\phi(\alpha), \sigma) \leq L_n\right\}\\
&\ge\#\left\{\phi\in\Map(S)\text{ with } \phi^{-1}(\sigma)\in\calD_\alpha^\sigma\text{ and } \iota(\alpha, \phi^{-1}(\sigma)) \leq L_n\right\}\\
&=\#\left\{\phi\in\Map(S)\text{ with } \phi(\sigma)\in\calD_\alpha^\sigma\text{ and } \iota\left(\alpha, \frac{\phi(\sigma)}{L_n}\right) \leq 1\right\}\\
&=\left(\sum_{\phi\in\Map(S)}\delta_{\frac 1{L_n}\phi(\sigma)}\right)\left(\left\{\lambda\in\calD_\alpha^\sigma,\ \iota(\alpha, \lambda) \leq 1\right\}\right)\\
&=\FM_{\sigma}^{L_n}\left(\left\{\lambda\in\calD_\alpha^\sigma,\ \iota(\alpha, \lambda) \leq 1\right\}\right).
\end{align*}
In light of \eqref{eq1} we get with $\Lambda$ as in \eqref{eq111} that
$$\Lambda\ge\liminf_{n\to\infty}\FM_\sigma^{L_n}(\{\lambda\in\calD_\alpha^\sigma,\ \iota(\alpha, \lambda) \leq 1\})\ge c\cdot\FM_{\Thu}(\{\lambda\in\calD_\alpha^\sigma,\ \iota(\alpha, \lambda) \leq 1\})$$
An analogue computation using the right hand side of \eqref{eq2} implies
$$\Lambda\le\limsup_{n\to\infty}\FM_\sigma^{L_n}(\{\lambda\in\bar\calD_\alpha^\sigma,\ \iota(\alpha, \lambda) \leq 1\})\le c\cdot\FM_{\Thu}(\{\lambda\in\bar\calD_\alpha^\sigma,\ \iota(\alpha, \lambda) \leq 1\})$$
Now, \eqref{eq000} implies that $\FM_{\Thu}(\bar\calD_\alpha^\sigma\setminus\calD_\alpha^\sigma)=0$ and hence that both bounds agree. We have proved that
$$\Lambda=c\cdot\FM_{\Thu}(\{\lambda\in\calD_\alpha^\sigma,\ \iota(\alpha, \lambda) \leq 1\})$$
and hence that
\begin{equation}\label{eq3}
c=\frac{\Lambda}{\FM_{\Thu}(\{\lambda\in\calD_\alpha^\sigma,\ \iota(\alpha, \lambda) \leq 1\})}.
\end{equation}
As we wanted, we have proved that $c$ in \eqref{eq1} does not depend on the sequence $(L_n)$ in question. In other words we have proved that
\begin{equation}\label{eq fifth turn}
\lim_{L\to\infty}\FM_\sigma^{L}=c\cdot\FM_{\Thu}
\end{equation}
exists for our particular filling multicurve $\sigma$. It remains to prove that it also exists for any other $\gamma$.

Now, repeating the same computations as above but with $\gamma$ instead of $\alpha$ we get that 
$$\FM_\sigma^L(\{\lambda\in\calD_{\gamma}^\sigma,\ \iota(\gamma,\lambda)\le 1\})\le \#\{\gamma'\text{ of type }\gamma\text{ with }\iota(\gamma',\sigma)\le L\}\le\FM_\sigma^L(\{\lambda\in\bar\calD_{\gamma}^\sigma,\ \iota(\gamma,\lambda)\le 1\})$$
Taking \eqref{eq fifth turn} into account we can pass to limits and get
\begin{align*}
c\cdot\FM_{\Thu}(\{\lambda\in\calD_{\gamma}^\sigma,\ \iota(\gamma,\lambda)\le 1\})&\le\liminf_{L\to\infty}\frac 1{L^{6g-6}} \#\{\gamma'\text{ of type }\gamma\text{ with }\iota(\gamma',\sigma)\le L\}\\
\limsup_{L\to\infty}\frac 1{L^{6g-6}}\#\{\gamma'\text{ of type }\gamma\text{ with }\iota(\gamma',\sigma)\le L\}&\le c\cdot\FM_{\Thu}(\{\lambda\in\bar\calD_{\gamma}^\sigma,\ \iota(\gamma,\lambda)\le 1\})
\end{align*}
Now, using again that the $\FM_{\Thu}(\bar\calD_\gamma^\sigma\setminus\calD_\gamma^\sigma)=0$ we get that the $\liminf$ and $\limsup$ agree and hence that 
$$\lim_{L\to\infty}\frac 1{L^{6g-6}}\#\{\gamma'\text{ of type }\gamma\text{ with }\iota(\gamma',\sigma)\le L\}= c\cdot\FM_{\Thu}(\{\lambda\in\calD_{\gamma}^\sigma,\ \iota(\gamma,\lambda)\le 1\})$$
We have thus proved that the desired limit exits for the positive continuous homogenous function
$$f_0:\CC(S)\to\BR_+,\ \ f_0(\lambda)=\iota(\lambda,\sigma)$$
We thus get from Proposition \ref{prop2} that the limit
$$\lim_{L\to\infty}\frac 1{L^{6g-6}}\#\{\gamma'\text{ of type }\gamma\text{ with }f(\gamma')\le L\}$$
exists for any other positive continuous homogenous function $f:\CC(S)\to\BR_+$. It exists in particular for the hyperbolic length fuction $\ell_X$, proving Theorem \ref{Mirzakhani-general}.\qed

\section{What changes without our simplifying assumptions?}
So far we assumed that the surface $S$ was closed and of genus $g\ge 3$. Let us first see what would change if the surface had cusps. In this case we would have to work on some compact surface $\Sigma$ with boundary and with interior homeomorphic to $S$. The reason why we would have to replace $S$ by $\Sigma$ is that quite a few things go wrong when we work with currents on non-compact surfaces. For example, continuity of the the intersection form fails when we work with currents on non-compact finite area hyperbolic surfaces \cite{Marie}. Such problems disappear if one works with currents on the compact surface $\Sigma$. There are however more annoying little issues waiting for us. For example, one needs to make sure that all the measures we are getting as limits are not supported by the boundary components of $\Sigma$. This is why in \cite{book}, where Theorem \ref{Mirzakhani-general} is proved in all generality, we actually worked with measures supported by $\CC_K(S)$, the set of current whose support projects into some compact set $K\subset\Sigma\setminus\D\Sigma$. Once these technical issues are dealt with, the argument in Section \ref{sec proof} applies word-by-word...

...at least as long as the surface is not exceptional, that is, as long as the center $\mathrm{Ker}$ of the mapping class group is trivial---this was guaranteed here by the assumption that $S$ has at least genus $3$. If $\mathrm{Ker}$ is not trivial then one should first assume that the multicurves under consideration are $\mathrm{Ker}$-invariant. Under this assumption, one can work on the action of $\Map/\mathrm{Ker}$ on the space $\CC^{\mathrm{Ker}}$ of $\mathrm{Ker}$-invariant currents, and the argument remains identical to the one we presented above. To deal with the case that the involved curves are not $\mathrm{Ker}$-invariant one can then use the same argument as in \cite[p.146]{book}. To conclude, the notation $\mathrm{Ker}$ is explained by the fact that the center of the mapping class group agrees with the kernel of its action on $\CM\CL(S)$.

\bibliographystyle{abbrv}

\end{document}